\font\titlefont= cmcsc10 at 12pt
\def\F{{\bf F}}
\def\P{{\bf P}}
\font\eighteenbf=cmbx10 scaled\magstep3
\def\vandaag{\number\day\space\ifcase\month\or
 januari\or februari\or  maart\or  april\or mei\or juni\or  juli\or
 augustus\or  september\or  oktober\or november\or  december\or\fi,
\number\year}
\magnification\magstep1
\vglue 2.0cm
\centerline{\eighteenbf Kummer Covers with Many Points}
%\bigskip
%\centerline{\eighteenbf }
%\bigskip
%\centerline{\eighteenbf  Titel}
\bigskip
\vskip 2pc
\centerline{\titlefont Gerard van der Geer \& Marcel van der Vlugt
} %%\footnote{}{\tt kummer, \vandaag \space (\number\time) }
\bigskip\bigskip
\centerline{\bf  Introduction}
\bigskip
\noindent
Let $p$ be a prime and $\F_q$ be a finite field with $q=p^m$ elements and let
$\bar{\F}_q$ be an algebraic closure of $\F_q$. In this paper we present a
method for constructing curves over finite fields with many points which are
Kummer covers of $\P^1$, or of other suitable base curves. For this we
look at rational functions $f \in {\F}_q(x)$ such that $f$ assumes a fixed
value $a \in \F_q$ on a (preferably large) set ${\cal P} \subseteq \P^1(\F_q)$.
To take a concrete example we set $a=1$. Then the algebraic curve which is the
Kummer covering  of $\P^1$ given by the equation 
$$
y^{q-1}=f(x),
$$
has fibres with many rational points and judicious choices of $f$ lead to
improvements and extensions of the tables [2] of curves with many points.
The methods we employed in the past were mostly based on Artin-Schreier covers
of the projective line.

In Section 1 we sketch the method and describe a way to obtain good rational
functions as above. This is based on an appropriate splitting
$f=f_1+f_2$ of a linearized polynomial $f$ having as zero
set a linear subspace $L$ of $\F_q$. In the following section we treat the case
where the linear space $L$ is the full space $\F_q$.  We find curves
$C_m$ defined over $\F_q$ for which the ratio $\# C_m(\F_q)/g(C_m)$  of the
number of rational points by the genus exceeds $\sqrt{q}$ for $m$ even and 
$2\sqrt{pq}/(p+1)$ for $m$ odd.  For $g$
large compared to $q$  the only way known so far to prove the existence of
curves with a comparable ratio is by class field theory, which is less
explicit (cf.\ [1]). Note that the result of
Drinfeld-Vladut, that $\lim \sup_{g \to \infty} \#C(\F_q)/g(C)\leq
\sqrt{q}-1$, shows that for given $q$ there are only finitely many
isomorphism classes of curves $C$ of  over $\F_q$  whose ratio
$\#C(\F_q)/g(C)$ exceeds $\sqrt{q}$ for $m$ even.

In  section 3 we consider the case where the linear subspace is of
codimension $1$ in $\F_q$ and we complement this note with a section with
variations.  We illustrate the sections with numerous examples and  thus obtain
a number of improvements of the existing tables. In many cases the methods also
give a  relatively easy way to construct for  certain pairs $(q,g)$ a curve
realizing the lower entry of the interval in the tables [2].  We conclude the
article with  two tables summarizing the new results from the examples. 

\vfill\eject
\bigskip
\centerline{ \bf \S 1. The Method }
\smallskip
\noindent
We consider the non-singular complete irreducible Kummer curve $C$ over $\F_q$
defined by the affine equation
$$
y^{q-1}= f(x),
$$
where the rational function $f(x) \in F_q(x)$ satisfies the
following conditions.

\smallskip
\proclaim (1.1) Conditions. 
\item{{\bf i)}} $f$ is not the $d$-th power of an element $g \in \bar{\F}_q(x)$ for
any divisor $d>1$ of $q-1$;
\item{{\bf ii)}} $f(x)=1$ on a substantial subset ${\cal P}$
of $\P^1(\F_q)$;
\item{{ \bf iii)}} $f(x)$ has many multiple zeros and poles.
\par
\smallskip
By i) the curve $C$ is a cyclic cover of $\P^1$ of degree $q-1$, by ii) the
curve $C$ has at least $(q-1)\# {\cal P}$ rational points and  condition iii)
keeps the genus of $C$ within bounds. 

The Hurwitz-Zeuthen formula gives the genus of $C$ (cf.\ e.g.\ [3]):

\proclaim (1.2) Proposition. If the divisor of
$f$ is  $(f)= \sum_{i=1}^\ell d_i P_i$ with distinct $P_i \in \P^1(\bar{\F}_q)$
then the genus $g(C)$ of $C$ is given by
$$
2g(C)-2= (\ell-2)(q-1) -\sum_{i=1}^\ell \gcd (q-1, |d_i|). \eqno(1)
$$
\par
\noindent
Note that a small value of $\ell$ and the greatest common divisors influence
the genus in a favourable way for our game.
\bigskip
Rational functions which satisfy Conditions (1.1)  arise for
instance in the following way. 

Let $L$ be an $r$-dimensional subspace of the $\F_p$-vector space
$\F_{q=p^m}$ with $r\geq 2$. Then the polynomial
$$
R=\prod_{c \in L} (x-c)
$$
is a $p$-linearized polynomial, i.e.\ of the form
$$
R(x)= \sum_{i=0}^r a_ix^{p^i} \quad \in \F_q[x]
$$
and moreover satisfies $a_0a_r\neq 0$. 

Now we split $R$ as
$$
R(x)=R_1(x)+R_2(x)\eqno(2)
$$
such that $R_1(x)= \sum_{i=s}^r b_i x^{p^i} \in \F_q[x]$ and $R_2(x)=
\sum_{i=0}^t c_i x^{p^i}$ with $0<s<r$, $t\leq s$, $b_sb_r\neq 0$ and $c_0c_t
\neq 0$. We denote the zero sets of $R_1$ (resp.\  $R_2$) by $L_1$ (resp.\
$L_2$) with $\# L_1= p^{r-s}$ (resp.\ $\#L_2= p^t$). Furthermore, in
connection with Condition (1.1) i) we require that $L_1 \neq L_2$. 

It is  obvious that
$$
f(x)= - {R_1(x)\over R_2(x)} = - {\big(\sum_{i=s}^r b_i^{1/{p^s}}
x^{p^{i-s}}\big)^{p^s}
\over \sum_{i=0}^t c_i x^{p^i}} \eqno(3)
$$
satisfies $f(x)=1$ for $x \in L-(L_1\cup L_2)$.  From (2) it follows that
$L\cap (L_1\cup L_2) = L_1\cap L_2$, which means that $L-(L_1\cup L_2)= L-
(L_1\cap L_2)$. Moreover, the zeros of $R_1$ and the pole $\infty$ have
multiplicities $>1$. Hence $f$ satisfies the Conditions (1.1).

\proclaim (1.3) Proposition. The Kummer cover $C$ of $\P^1$ defined by
the equation $y^{q-1}=f(x)$ with $f(x)= -R_1/R_2$ as in (3) has genus
$$
g= \big\{ (p^{r-s}+p^t-\delta-1)(q-2) -\delta p^{{\rm
gcd}(m,s)}-p^{{\rm gcd}(m,r-t)}+2\delta +2 \big\}/2 \eqno(4)
$$
and the number of $\F_q$-rational points on $C$ satisfies
$$
\# C(\F_q)\geq (p^r-\delta)(q-1), \eqno(5)
$$
where $\delta= \# (L_1\cap L_2)$. \par
\smallskip
\noindent
{\sl Proof.} By the assumption $L_1\neq L_2$ it follows that the function $f$
satisfies (1.1) i). The divisor of $f$ is
$$
(f)= \sum_{P \in L_1\cap L_2} (p^s-1)P + \sum_{P \in L_1-(L_1\cap L_2)} p^s P
- \sum_{P \in L_2-(L_1\cap L_2)} P - (p^r-p^t)P_{\infty}.
$$
The number $\ell$ of distinct zeros and poles of $f$ is 
$$
\# L_1 + \# L_2 - \# (L_1 \cap L_2)+1=
p^{r-s}+p^t -\delta +1.
$$
According to Prop.\ (1.2) the genus satisfies
$$
2g(C)-2= (p^{r-s}+p^t-\delta -1)(q-1)-\delta(p^{{\rm
gcd}(m,s)}-1)-(p^{r-s}-\delta)- (p^{{\rm gcd}(m, r-t)}-1)
$$
and we obtain (4). For $x \in L-(L_1\cap L_2)$ we have $f(x)=1$ and  thus
 over each $y\in \F_q^*$ we find $p^r-\delta$ rational points on $C$. Other
rational points could come from the branch points of $C$. The set of branch
points is $L_1 \cup L_2 \cup \infty$ and they contribute  rational points if
the ramification points over such branch points happen to be rational. This
yields the required estimate (5). $\bullet$
\smallskip
\noindent
{\bf (1.4) Example.} Take
$\F_{16}$ with
$L=\F_{16}$. Then
$R=x^{16}+x$ and we split
$R$ as $R=R_1+R_2$ with $R_1=x^{16}+x^2$ and $R_2= x^2+x$. In this case $r=4$,
$s=t=1$, $L_1=\F_8$, $L_2=\F_2$ and $\delta=2$. From Prop.\ (1.3) we see that
the curve $C$ defined over $\F_{16}$ by 
$$
y^{15}=(x^{16}+x^2)/(x^2+x)= x^{14}+x^{13}+\ldots +x
$$
has genus $g(C)= 49$ and $\# C(\F_{16})= 14\times 15 + 3=213$ since the
ramification points over the branch points in $\F_2 \cup \infty$ are rational.
This provides a new entry for the tables in [2].
\bigskip
We remark that the ratio $\#C(\F_q)/g(\F_q)$ for the curves that appear in
Proposition (1.3) exceeds $2p^r/(p^{r-s}+p^t)$,  which is
optimal for $s=t=[r/2]$. For that choice
$$
\# C(\F_q)/g(C) > \cases{ \sqrt{p^r}& for $r$ even, \cr
 2\sqrt{p^{r+1}}/(p+1)& for $r$ odd.\cr} \eqno(6)
$$
From (6) it follows that the case $L=\F_{p^m}$ with $R= x^{p^m}-x$ is of
special interest.

\bigskip
\centerline{\bf \S 2. The Case $L=\F_q$}
\bigskip
\noindent
In this Section we consider the case where $L$ equals the full vector space
$\F_q$.
For odd $m$ we write
$$
x^{p^m}-x= R_1+R_2= (x^{p^m}-ax^{p^{(m-1)/2}})+(ax^{p^{(m-1)/2}}-x),
$$
with $a \in \F_q^*$, i.e., we look at the case $s=t=[m/2]$. Since
$$
{\rm gcd}(x^{p^m}-ax^{p^{(m-1)/2}}, ax^{p^{(m-1)/2}}-x)={\rm gcd}(x^{p^m}-x,
ax^{p^{(m-1)/2}}-x)
$$ 
we have for $u \in \F_q^*$
$$
u \in L_1\cap L_2 \iff u^{p^{(m-1)/2}-1}= 1/a.
$$
This equation has no solutions in $\F_q^*$ if $a$ is not a 
$(p^{(m-1)/2}-1)$-th power in $\F_q^*$ and the number of solutions in $\F_q^*$
is ${\rm gcd}(p^{(m-1)/2}-1, p^m-1)= p-1$ if $a$ is a $(p^{(m-1)/2}-1)$-th
power in $\F_q^*$. The latter holds always if $p=2$.

First we consider the case that $a$ is a $(p^{(m-1)/2}-1)$-th
power in $\F_q^*$. Often we shall write $a \in (\F_q^*)^d$ to indicate that $a$
is a $d$-th power in $\F_q^*$.

\proclaim (2.1) Proposition. For odd $m\geq 3$ the curve $C_m$ defined over
$\F_{q=p^m}$ by the equation 
$$
y^{q-1}= - { (x^{p^{(m+1)/2}}-a^{p^{(m+1)/2}}x)^{p^{(m-1)/2}} \over
ax^{p^{(m-1)/2}}-x} 
$$
with $a \in (\F_q^*)^{p^{(m-1)/2}-1}$
has genus 
$$
g(C_m)= \big\{ (p^{(m+1)/2} + p^{(m-1)/2} -p-1)(q-2) -p^2+p+2\big\} /2 
$$
and has the following number of rational points
$$
\# C_m(\F_q)= \cases{(q-1)(q-p) & for odd $p$,\cr
(q-1)(q-p)+3 & for $p=2$. \cr} 
$$\par
\smallskip
\noindent
{\sl Proof.} The degree of ${\rm gcd}(R_1,R_2)={\rm
gcd}(x^{p^{(m+1)/2}}-a^{p^{(m+1)/2}}x, ax^{p^{(m-1)/2}}-x)$ is the cardinality
of the $\F_p$-vector space $L_1\cap L_2$. The condition
that $a$ is $p^{(m-1)/2}-1$-th power implies $\delta=p$.
We have  $s=t=(m-1)/2$ and the expression
for the genus follows now directly by substitution in Prop.\ (1.3). Over each
$y  \in \F_q^*$ we have $p^m-\delta=p^m-p$ rational points on $C_m$. The only
branch points which possibly contribute rational points to $C_m$
are the branch points in $\F_p \cup \infty$.  Over each point of
$\F_p\cup \infty$ there lie $p-1$ ramification points on $C_m$. These are
rational if and only if $-a$ is a $(p-1)$-th power in $\F_{q}$. This holds for
pairs $(p,m)$ with $pm$ even  which implies our formula for $\# C_m(\F_q)$.
$\bullet$

\bigskip
\noindent
{\bf (2.2)  Example.} As an illustration  of  Proposition (2.1) we take $p=3$,
$m=3$ and get the curve $C$ over $\F_{27}$  given by
$$-y^{26}= x^{24}+x^{22}+\ldots + x^2
$$
with $g(C)=98$ and $\# C(\F_{27})=624$. In this case the Oesterl\'e upper bound
is
$b=745$, so $C$ satisfies our qualification criterion $\#C(\F_{27}) \geq
[b/\sqrt{2}]$ for the tables in [2].

For another example we take $\F_{32}$. Then the curve $C$ with affine equation
$$
y^{31}= (x^8+x)^4/(x^4+x)
$$
has genus $g(C)= 135$ and $\# C(\F_{32})= 31 \times 30 + 3 =933$. The
Oesterl\'e upper bound in this case is $1098$. 

For $q=3^5$ we obtain from Prop.\ (2.1) a curve $C$ of genus $g(C)=3854$ and
$\# C(\F_{243})=58080$. The Oesterl\'e upper bound is $81835$.
\smallskip

\bigskip
For $a$ not a $(p^{(m-1)/2}-1)$-th power in $\F_q^*$ we have a similar
proposition.
\proclaim (2.3) Proposition. For odd $m\geq 3$ the curve $C_m$ over
$\F_{q=p^m}$ defined by
$$
y^{q-1}= -{ (x^{p^{(m+1)/2}}-a^{p^{(m+1)/2}}x)^{p^{(m-1)/2}} \over
ax^{p^{(m-1)/2}}-x}
$$
with $a\not\in (\F_q^*)^{p^{(m-1)/2}-1}$ has genus
$$
g(C_m)=\big\{ (p^{(m+1)/2}+p^{(m-1)/2}-2)(q-2)-2p+4\big\} /2
$$
and has the following number of rational points
$$
\# C_m(\F_q)= \cases{(q-1)^2& if $-a \not\in (\F_q^*)^{p-1}$,\cr
(q-1)^2+2(p-1)& if $-a \in (\F_q^*)^{p-1}$.\cr}
$$
\par
\smallskip
\noindent
{\sl Proof.} The proof is similar to that of Prop.\ (2.1). with the following
modifications. In this case ${\rm gcd}(R_1,R_2)$ has degree $1$ which means
that $\delta= \# (L_1 \cap L_2)=1$ and over each $y\in \F_q^*$ we have
$p^m-\delta=p^m-1$ rational points on $C_m$. The branch points which possibly
contribute rational points on $C_m$ are $0$ and $\infty$.  Over these points
there are $p-1$ ramification points on $C_m$ which are rational points if and
only if $-a$ is a $(p-1)$-th power in $\F_q$. This gives the formula for the
number of rational points. $\bullet$

\smallskip
\noindent
{\bf (2.4) Examples.}  For $p$ odd  we take $a=-1$ since $-1$ is not a
$(p^{(m-1)/2}-1)$-th power in $\F_q^*$. Over $\F_{27}$ the curve $C_m$ has
genus $g(C_m)= 124$ and
$\#C_m(\F_{27})=680$ while the  Oesterl\'e upper bound is $901$.
Over $\F_{3^5}$ we find $g(C_m)=4096$ and $\#C_m(\F_{3^5})=58568$. The
Oesterl\'e upper bound is here $86441$.

\bigskip

For $q=p^m$ with $m$ even the splitting
$$
x^q-x= (x^q-ax^{\sqrt{q}})+(ax^{\sqrt{q}}-x),
$$
where $a \in \F_q^*$ is such that $a\notin (\F_q^*)^{\sqrt{q}-1}$ yields very
good curves.

\proclaim (2.5) Proposition. If $q=p^m$ with $m$ even then the
curve
$C_m$ defined over $\F_{q=p^m}$ by the equation
$$
y^{q-1}= - { x^q-ax^{\sqrt{q}} \over ax^{\sqrt{q}} -x}
\qquad {\rm with} \quad a 
\in \F_q^*, a\not \in (\F_q^*)^{\sqrt{q}-1} 
$$
has genus $g(C_m)= (\sqrt{q}-1)(q-2) - \sqrt{q}+2$ and $\#C_m(\F_q)= (q-1)^2$.
\par
\smallskip
\noindent
{\sl Proof.} In this situation we have $s=t=m/2$ and the condition
$a^{\sqrt{q}+1}\neq 1$ implies $L_1\cap L_2 = \{ 0 \}$, so $\delta=1$. The
formula for $g(C_m)$ follows from Prop.\ (1.3). Over each $y \in \F_q^*$ there
are $p^m-\delta$ rational points on $C_m$.  The only branch points which
possibly give rise to rational points on
$C_m$ are
$0$ and
$\infty$. The ramification points over $0$ (resp.\ $\infty$) on $C_m$ are
rational iff the equation $w^{\sqrt{q}-1}=-a$  (resp.\ $w^{\sqrt{q}-1}=
(-1/a)$) is solvable in $\F_q$. Since  $a^{\sqrt{q}+1}\neq 1$
these equations have no solutions in $\F_q$  and consequently we have
$\#C_m(\F_q)= (q-1)^2$.
$\bullet$
\smallskip
Note that in this case 
$$
\# C_m(\F_q)/g(C_m) > \sqrt{q} +1.
$$
\smallskip
\noindent
{\bf (2.6) Examples.} For $q=9$ we find $g(C_m)= 13$ and $\#C_m(\F_9)= 64$.
This is very close to the Oesterl\'e upper bound $66$ and might well be optimal
(i.e.\ equal to the actual maximum number $N_q(g)$, cf.\ [2]).
For $q=16$ we find $g(C_m)=40$, $\#C_m(\F_{16})=225$; the Oesterl\'e upper
bound is $244$.  For $q=64$ we find
$g(C_m)= 428$ and $\#C_m(\F_{64})=3969$ with Oesterl\'e upper bound $4786$.
For $q=81$ we find $g(C_m)=625$ and $\#C_m(\F_{81})=6400$, still reasonable
compared with the Oesterl\'e upper bound $7824$.
\bigskip
\centerline{\bf \S 3. Subspaces of  Codimension $1$ }
\bigskip
\noindent
We take as subspace of $\F_q$  the $(m-1)$-dimensional
subspace
$$
L= \{ x \in \F_q: {\rm Tr}_{\F_q/\F_p}(x)=0\},\qquad {\rm where} \quad 
{\rm Tr}_{\F_q / \F_p}(x)= x^{p^{m-1}}+\ldots + x^p + x,
$$
and put $R(x)=\sum_{i=0}^{m-1} x^{p^i}$. Note that by a transformation $x
\mapsto ax$ on $\F_q$ with $a \in \F_q^*$ we can transform any codimension $1$
space into this subspace $L$. We split the polynomial $R$ as $R_1+R_2$  with
$R_1=
\sum_{i=s}^{m-1} x^{p^i}$ and $R_2=\sum_{i=0}^{s-1} x^{p^i}$. The corresponding
curve $C_m$ over $\F_q$ is defined by
$$
y^{q-1}= -(x^{p^{m-1-s}}+\ldots +x)^{p^s}/(x^{p^{s-1}}+\ldots + x) \eqno(7)
$$
Applying Prop.\ (1.3) to this situation gives the following result.

\proclaim (3.1) Proposition. For $m\geq 3$ and $0< s < m-1$ such that ${\rm
gcd}(m,s)=1$ the curve $C_m$ given by (7) has genus
$$
g(C_m)= \big\{ (p^{m-1-s}+p^{s-1} -2)(q-2) -2p+4\big\} /2
$$
and
$$\# C_m(\F_q)=\cases{(p^{m-1}-1)(q-1) & if $pm$ odd and $p
\not| s(m-s)$,\cr 
(p^{m-1}-1)(q-1)+(p-1) & if $pm$ odd  and $p| s(m-s)$,\cr
(p^{m-1}-1)(q-1)+2(p-1) & if $pm$ even and $p \not|
s(m-s)$,\cr 
(p^{m-1}-1)(q-1)+3(p-1) & if $pm$ even and $p|
s(m-s)$. \cr}
$$
\par
\smallskip
\noindent
{\sl Proof.} In the notation used in Section 1 we find 
$$
L_1= \{ x \in \F_{p^{m-s}} :
{\rm Tr}_{\F_{p^{m-s}}/\F_p} (x)=0 \}
$$
and
$$
L_2= \{ x \in \F_{p^{s}} :
{\rm Tr}_{\F_{p^{s}}/ \F_p} (x)=0 \}.
$$
Then $L_1 \cap L_2\subset \F_{p^{m-s}}\cap \F_{p^s}=\F_{p^{{\rm
gcd}(m,s)}}=\F_p$. Combining ${\rm gcd}(m,s)=1$ with the condition on the
traces gives $L_1 \cap L_2=\{ 0 \}$, hence $\delta =1$. If $p\not|
s(m-s)$  the ramification points over $L_1-\{ 0 \}$ and $L_2-\{0 \}$ are not
rational. On the other hand, if $p| s(m-s)$ then $\F_q^* \cap (L_1 \cup L_2)=
\F_p^*$ and the ramification points over $\F_p^*$ are rational.  Over the
branch point
$0$ (resp.\ 
$\infty$)  which has  multiplicity $p^s-1$ (resp. $p^{m-1}-p^{s-1}$) there lie
$p-1$ ramification points on $C_m$. These are rational if and only if $-1$ is
a $(p-1)$-th power in $\F_q$ which holds if and only if $pm$ is even. The
formulas now follow from Prop.\ (1.3). $\bullet$
\bigskip

\bigskip
\noindent
{\bf (3.2) Examples.} Take $\F_{27}$ then $L=\{ x \in \F_{27}: {\rm
Tr}_{\F_{27}/\F_{3}}(x)=0 \}$ is given by $R(x)=x^9+x^3+x$ and we can consider
the curve
$$
C: \quad y^{26}= - (x^8+x^2).
$$
It follows from Prop.\ (3.1) that $g(C)=24$ and $\#C(\F_{27})=208$ which
improves [2].

For $\F_{32}$ with (7) of the form $y^{31}=(x^4+x^2+x)^4/(x^2+x)$ we obtain
according to Prop.\ (3.1) a curve $C$ of genus $60$ and $\# C(\F_{32})=468$.
The Oesterl\'e upper bound is $542$.

\smallskip
\noindent{\bf (3.3) Example.} Finally we consider an example where ${\rm
gcd}(m,s)\neq 1$. Take
$\F_{64}$ and
$$
f(x)= {x^{32}+x^{16} \over x^8+x^4+x^2+x} = {(x^2+x)^{15}\over
(x^4+x+1)(x^2+x+1)}.
$$
For the curve $C$ given by $y^{63}=f(x)$ Prop.\ (1.2) implies that
$2g-2=7\times 63 - 3 \times 3 - 6 \times 1$, hence $g(C)=214$. Each of the
branch points $0$,
$1$ and $\infty$ induces $3$ rational ramification points on $C$, the zeros of
$x^2+x+1$ induce $1$ ramification point each, while the ramification points
from the zeros of $x^4+x+1$ are not rational. The number of rational points on
$C$ is thus $\# C(\F_{64})= (32-2)\times 63 + 11= 1901$. The Oesterl\'e upper
bound is $2553$.

\smallskip
\noindent
{\bf (3.4) Remark.} For even $m$ and $L=\{ x \in \F_{q=p^m} : {\rm
Tr}_{\F_{p^m}/\F_p} (x)=0\}$ the splitting
$$
\sum_{i=0}^{m-1} x^{p^{i}} = R_1+R_2= \sum_{i=m/2}^{m-1} x^{p^i} +
\sum_{i=0}^{(m/2)-1} x^{p^i}
$$
does not satisfy condition i). The corresponding equation $y^{q-1}= -R_1/R_2=
-R_2^{\sqrt{q}-1}$ leads to the curve
$$
C: \quad y^{\sqrt{q}+1}= aR_2= a(x^{p^{(m/2)}-1}+ x^{p^{(m/2)}-2}+\ldots + x)
\eqno(8)
$$
where $a \in \F_q^*$ is such that $a^{\sqrt{q}}+a=0$.

To determine $g(C)$ and $\#C(\F_q)$ we consider the $\F_p$-linear map $\phi$
on $L$ defined by $\phi(x)=aR_2(x)$. The kernel of $\phi$ is
$$
\ker (\phi)= \{ x \in \F_{\sqrt{q}} : {\rm Tr}_{\F_{\sqrt{q}}/\F_p}(x)=0\}
\quad {\rm and} \quad \phi(L)= \F_{\sqrt{q}}.
$$
For $y \in \F_q^*$ we have $y^{\sqrt{q}+1}\in \F_{\sqrt{q}}^*$ so over each $y
\in \F_q^*$ there are $\# \ker (\phi) = \sqrt{q}/p$ rational points on $C$.
The set of branch points is $\ker(\phi)\cup \infty$ and each branch point
induces $1$ rational point on $C$. Hence
$$
\# C(\F_q)= (q-1)\sqrt{q}/p + \sqrt{q}/p +1 = (q \sqrt{q}/p )+1.
$$
From Prop.\ (1.2) we find $g(C)= (q-p\sqrt{q})/2p$. We thus get explicit
maximal curves:
\bigskip
\proclaim (3.5) Proposition. The curve $C$ over $\F_q$ given by (8) with
$g(C)= (q-p\sqrt{q})/2p$ and $\#C(\F_q)= (q\sqrt{q}/p)+1$ is a maximal curve,
i.e.\ it attains the Hasse-Weil upper bound. \par

By the substitution  $x \mapsto z^p-z$ in (8) we obtain the equation for the
Hermitian curve $y^{\sqrt{q}+1}=a(z^{\sqrt{q}}-z)$. So the curve $C$ figuring
in Prop.\ (3.5) is a  quotient of the Hermitian curve.

\bigskip
\centerline{\bf \S 4. Variations}
\bigskip
To find curves with many points with this method it is not necessary to depart
from a linearized polynomial. This is illustrated by the following example
where we take a curve of the form
$$
y^{q-1}= xf(x)^p
$$
with $f(x)\in \F_q[x]$.
\smallskip
\noindent
{\bf (4.1) Example}. Take $\F_{16}$ and consider the irreducible complete
non-singular curve $C$ given by the affine equation
$$
y^{15}= x(x^2+x+1)^2.
$$
Remark that $x(x^2+x+1)^2=x^5+x^3+x$ satisfies  Conditions (1.1).
According to (1) the curve $C$ has genus $g(C)=12$ and the number of points is
$\# C(\F_{16})= 15\times 5 + 8=83$, where the branch point 
$\infty$ contributes $5$ rational points and the branch points in
$\F_{4}-\{1\}$ each contribute $1$ rational point. This example provides a new
entry for the tables [2], where the interval [68--97] is given.

\bigskip
An advantage of our method is that we can also find  good curves $C$ such that
only a few fibres over $\P^1(\F_q)$ contribute to the rational points on $C$,
but these then do so substantially, as in the preceding example.  We can use
this for instance to construct Artin-Schreier covers of
$C$ given by
$$
z^p-z= h(x),
$$
where in order to obtain good curves one has to impose the condition
${\rm Tr}(h(x))=0$ for a few values $x$ only.
\smallskip
\noindent
{\bf (4.2) Example}. Take the field $\F_{32}$ and consider the curve $C$
defined by
$$
y^{31}= x^5+x^3.
$$
The polynomial $x^5+x^3+1$ is irreducible over $\F_2$, so it has $5$ zeros
in $\F_{32}$. There are three ramification points $P_0,P_1$ and $P_{\infty}$ 
lying over
$0$,
$1$ and
$\infty$.

We find
$$
g(C)=15, \qquad \# C(\F_{32})=158,
$$
which comes up to the best value known for $(q,g)=(32,15)$ in [2].

We immediately see that the zeros $x \in \F_{32}$ of $x^5+x^3+1$ satisfy ${\rm
Tr}(x)=0$. The divisor of $x$ is
$$
(x)= 31 P_{0} - 31 P_{\infty}.
$$
The Artin-Schreier cover $\tilde{C}$ of $C$ given by
$$
z^2+z=x
$$
has $2$ rational points over each of the $155$ points $(x,y)$ of $C(\F_{32})$
with
$y\in \F_{32}^*$. We thus find
$$
\# \tilde{C}(\F_{32})= 2\times 155 + 1 +2=313,
$$
and $g(\tilde{C})= 45$. (See [3] for formulas for the
genus.) This improves [2], where the interval is [302--428].

As a variation on this theme we take $\F_{16}$ with the curve $C$ given by
$$
y^{15}= x^4+x^3.
$$
This has genus $g(C)=6$ with $65$ rational points. The Artin-Schreier cover
$\tilde{C}$ of $C$ defined by  $z^2+z= 1/x$ yields a curve of genus
$g(\tilde{C})=20$ with $\# \tilde{C}(\F_{16})=127$.

\bigskip
If one has a curve $C$ with many points then often a curve $C^{\prime}$
obtained as the image under a $\F_q$-morphism $C\to C^{\prime}$ is also a
good curve because the set of eigenvalues of Frobenius for $C^{\prime}$
is a subset of those for $C$. In the cases dealt with in the preceding sections
where the curve is of the form
$$
y^{q-1} = f(x^{p-1}),
$$
we can consider the curves $y^s= f(x^t)$ for any divisor $s$ of $q-1$ and $t$
of $p-1$.

\smallskip
\noindent{\bf (4.3) Example.} From the curve $C$ over $\F_{27}$ given in
Example (2.2) we obtain the curve $C^{\prime }$
$$
 -y^{13}= x^{24}+x^{22}+\ldots + x^2 \quad {\rm with}\quad
g(C^{\prime })=48\quad {\rm and} \quad \# C^{\prime }(\F_{27})= 316,
$$
where the tables give [325--402], and
$$
-y^{26}= x^{12}+x^{11}+\ldots + x \quad {\rm with } \quad  g(C^{\prime })=49
\quad{\rm and}
\quad \# C^{\prime }(\F_{27})= 314,
$$
a new entry in the tables.

\bigskip
Of course, the methods can be varied in several ways. For example, one can
replace $y^{q-1}$ by  $y^t$ for $t$ a divisor of $q-1$
and take a function $f$ which assumes for many $x$ a $t$-th power in $F_q$. We
now give an example of this.
\smallskip
\noindent
{\bf (4.4) Example}.  Take $\F_{81}$ and consider the curve given by the
equation
$$
y^{10}= x^2+x.
$$
For $y\in \F_{81}$ we have $y^{10}\in \F_9$, so the equation $x^2+x=y^{10}$
has always  solutions $x\in \F_{81}$. The curve $C$ has genus
$g(C)=4$ and $\#C(\F_{81})=154$. Consider the double cover $\tilde{C}$ of $C$
given by
$$
z^2=x^2+x+2.
$$
Over each $(x,y)\in C(\F_{81})$ with $y\in \F_{81}^*$ the curve $\tilde{C}$
has rational points since $x^2+x+2 \in \F_9$. A computation of the genus and
the number of points yields
$$
g(\tilde{C})= 17, \qquad \# \tilde{C}(\F_{81})=288.
$$
This is a new entry for the tables [2].
\bigskip
We can also apply the methods to a base curve different from $\P^1$ as the
following examples show.
\smallskip
\noindent
{\bf (4.5) Example.} Take $\F_8$ and consider the curve $C$ of genus $1$
defined by
$$
y^2+y=x+{1 \over x} +1.
$$
It has $14$ rational points, namely the two ramification points $P_0$ and
$P_{\infty}$, and six pairs of points $P_{\zeta}, P_{\zeta}^{\prime}$, one
over each $7$-th root $\zeta \neq 1$ of $1$. Consider now the cover $\tilde{C}$
of $C$ defined by 
$$
z^7= x(x^6+1)/(x+1), 
$$
cf.\ Prop.\ (2.1). It has branch points $P_0$ and $P_{\infty}$ and  $P_x,
P_x^{\prime}$ for $x$ a third root of unity.  Then the genus $g(\tilde{C})$
satisfies
$2g(\tilde{C})-2=7\times 0+8\times 6=48$, hence $g(\tilde{C})=25$. The rational
points come from $12$ fibres of order $7$  over $P_{\zeta}$ and
$P_{\zeta}^{\prime}$, and from the two ramification points over $P_0$ and
$P_{\infty}$, giving $\#\tilde{C}(\F_8)=86$ which improves the entry [84--97]
of the tables.

\smallskip
\noindent
{\bf (4.6) Example}. Take $\F_8$ and consider the Klein curve $C$ of genus $3$
defined by $y^3+x^3y+x=0$. It has $24$ rational points. Consider then the
cover $\tilde{C}$  given by  $z^7= x(x^6+1)/(x+1)$. The branch points on $C$
are the points lying over
$x=0$,
$x=$ a third root of unity and $x=\infty$. We find $g(\tilde{C})=51$ and
$\# \tilde{C}(\F_8)=132$. The Oesterl\'e upper bound is  $173$.
\smallskip
\noindent
{\bf (4.7) Example}. Take $\F_{9}= \F_3[i]$ with $i^2=-1$ and consider the
curve
$C$ of  genus $1$ defined by
$$
y^2=x^3+x.
$$
It has $16$ rational points over $\F_9$, the $4$ ramification points
$P_0$, $P_{\infty}$, $P_i$ and $P_{-i}$, and six pairs $P_x,P^{\prime}_x$ for
$x\in \F_9-\{0,\pm i\}$. Take the function  $f=x/y$ with divisor
$(f)=P_0+P_{\infty}-P_i-P_{-i}$ and consider the cover $\tilde{C}$ of
$C$ defined by
$$
z^4=f^3+f.
$$
Observe that for $u \in \F_9^*$ the expression $u^3+u$ is a $4$-th power in
$\F_9^*$. One has 
$(f(f^2+1))= P_0+P_{\infty}+2P_1+2P_1^{\prime}-3P_i-3P_{-i}$. The curve
$\tilde{C}$ has genus
$g(\tilde{C})=9$ and has
$10\times 4 + 4 + 4=48$ rational points. Here the points
$P_0,P_{\infty},P_i,P_{-i}$ are branch points with total ramification, while
the branch points $P_1$ and
$P_1^{\prime}$ each contribute $2$ rational points. This comes up to the best
known curve and is very close to the Oesterl\'e upper bound $51$.
\smallskip
\noindent
{\bf (4.8) Example}. Take $\F_{9}= \F_3[i]$ with $i^2=-1$ and consider the
curve
$C$ of genus $2$ defined by
$$
z^2=x(x^4+x^2+2).
$$
It has $18$ rational points over $\F_9$. We denote them by $P_0$,
$P_{\infty}$ and by  $P_x$, $P_x^{\prime}$ in the fibre over $x$ for each
$x\in \F_9^*$. According to Prop.\ (2.5) the Kummer cover
$D$ of $\P^1$ defined by 
$$
y^8= -{(x^9-ax^3)\over (ax^3-x)} = - {x^2 (x^2-a^3)^3 \over a
x^2-1}\qquad \hbox {\rm with $a$  such that }
\quad a^2+a+2=0 
\eqno(9)
$$
has fibres consisting of $8$ rational points over each $x \in \F_9^*$.
We consider the curve $\tilde{C}$ which is the cover of $C$ defined by (9).
The branch points on $C$ are the four points $P_0$, $P_{\infty}$, 
$P_{\xi}, P_{\xi}^{\prime}$ with $\xi^2=a^3$ and the four points $P_{\eta},
P_{\eta}^{\prime}$ with $\eta^2= 1/a$. The divisor of the function $f$ given by
the right hand side of (9) is
$$
(f)=
4P_0+6P_{\xi}+6P_{-\xi}-P_{\eta}-P_{\eta}^{\prime}-P_{-\eta}-P_{-\eta}^{\prime}-12
P_{\infty}.
$$
By Hurwitz-Zeuthen the genus is $33$. Over each $x \in \F_9^*$ we find $16$
rational points on $\tilde{C}$ giving $\# {\tilde{C}}(\F_9)=128$, a
significant improvement of the entry [109--133] in the tables [2].

If we take here instead of the base curve $C$ the curve $C^{\prime}$ of genus
$2$ with $18$ rational points defined by
$$
z^2=x(x^4+x^3+x^2+x+1)
$$
then (9) defines a Kummer cover $\tilde{C^{\prime}}$ of $C^{\prime}$ of genus
$41$ with $128$ rational points.
\smallskip

By employing the methods
in a systematic way we expect more improvements and supplements to the tables
in [2].

\bigskip
\vfill\eject 
\smallskip
\centerline{\bf \S 5. Tables}
\smallskip
\noindent
We now summarize the new results from our examples
for  tables of curves with many points.
\bigskip
%%%%%%%%%%%%%
%%%%%%%%%%%%%
\font\tablefont=cmr8
\def\quad{\hskip 0.6em\relax}
\def\quod{\hskip 0.6em\relax}
\def\vhop{
    height2pt&\omit&&\omit&&\omit&&\omit&\cr}
\noindent{\bf Table $p=2$.}
$$
\vcenter{
\tablefont
\lineskip=1pt
\baselineskip=10pt
\lineskiplimit=0pt
\setbox\strutbox=\hbox{\vrule height .7\baselineskip
                                depth .3\baselineskip width0pt}%
\offinterlineskip
\hrule
\halign{&\vrule#&\strut\quod\hfil#\quad\cr
\vhop
&$q$&&$g(C)$&&new entry&&old entry&\cr
\vhop
\noalign{\hrule}
\vhop
&8&&25&&[86--97]&&[84-97]&\cr
%&8&&48&&[131--164]&&[126--164]&\cr
&8&&51&&[132--173]&&&\cr
&16&&12&&[83--97]&&[68-97]&\cr
&16&&20&&[127--140]&&[121--140]&\cr
&16&&40&&[225--244]&&[197--244]&\cr
&16&&49&&[213--286]&&&\cr
&32&&45&&[313--428]&&[304--428]&\cr
&32&&60&&[468--542]&&&\cr
&32&&135&&[933--1098]&&&\cr
&64&&214&&[1901--2553]&&&\cr
&64&&428&&[3969--4786]&&&\cr
\noalign{\hrule}
\vhop
}
\hrule
}
$$
%%%

%%%%%%%%%%%%%
\font\tablefont=cmr8
\def\quad{\hskip 0.6em\relax}
\def\quod{\hskip 0.6em\relax}
\def\vhop{
    height2pt&\omit&&\omit&&\omit&&\omit&\cr}
\noindent{\bf Table $p=3$.}
$$
\vcenter{
\tablefont
\lineskip=1pt
\baselineskip=10pt
\lineskiplimit=0pt
\setbox\strutbox=\hbox{\vrule height .7\baselineskip
                                depth .3\baselineskip width0pt}%
\offinterlineskip
\hrule
\halign{&\vrule#&\strut\quod\hfil#\quad\cr
\vhop
&$q$&&$g(C)$&&new entry&&old entry&\cr
\vhop
\noalign{\hrule}
\vhop
&9&&13&&[64--66]&&[60--66]&\cr
&9&&33&&[128--133]&&[109--133]&\cr
&9&&41&&[128--158]&&[119--158]&\cr
&27&&24&&[208--235]&&[190--235]&\cr
&27&&49&&[314--409]&&&\cr
&27&&98&&[624--745]&&&\cr
&27&&124&&[680--901]&&&\cr
&81&&17&&[288--387]&&&\cr
&81&&625&&[6400--7824]&&&\cr
&243&&3854&&[58080--81835]&&&\cr
&243&&4096&&[58568--86441]&&&\cr
\noalign{\hrule}
\vhop
}
\hrule
}
$$
%%%
\vfill\eject
\bigskip
\centerline{\bf References}
\bigskip
\noindent
\smallskip
\noindent
[1] R.\ Auer: Ray class fields of global function fields with many rational
places. Preprint University of Oldenburg, 1998.
\smallskip
\noindent
[2] G.\ van der Geer, M.\ van
der Vlugt: Tables of curves with many points. To appear in {\sl Math. Comp.\ }
(1999).
\smallskip
\noindent
[3] H.\ Stichtenoth: Algebraic function fields and codes.  Springer Verlag,
Berlin 1993.
\bigskip
%%%
%%%
\bigskip
\settabs3 \columns
\+G. van der Geer  &&M. van der Vlugt\cr
\+Faculteit
Wiskunde en Informatica &&Mathematisch Instituut\cr
\+Universiteit van
Amsterdam &&Rijksuniversiteit te Leiden \cr
\+Plantage Muidergracht 24&&Niels Bohrweg 1 \cr
\+1018 TV Amsterdam
&&2333 CA Leiden \cr
\+The Netherlands &&The Netherlands \cr
\+{\tt geer@wins.uva.nl} &&{\tt vlugt@wi.leidenuniv.nl} \cr
\bye